\documentclass[reqno]{amsart}

\newtheorem{Theorem}{Theorem}[section]

\newtheorem{Corollary}[Theorem]{Corollary}

\theoremstyle{remark}
\newtheorem{Remark}[Theorem]{Remark}
\numberwithin{equation}{section}

\allowdisplaybreaks

\begin{document}

\title[Generalizations of Jacobi's triple product identity]
  {Abel--Rothe type generalizations of\\
   Jacobi's triple product identity}
\author{Michael Schlosser}

\address{Institut f\"ur Mathematik der Universit\"at Wien,
Strudlhofgasse 4, A-1090 Wien, Austria}
\email{schlosse@ap.univie.ac.at}
\urladdr{http://www.mat.univie.ac.at/{\textasciitilde}schlosse}
\thanks{The author was supported by an APART grant of the Austrian
Academy of Sciences}
\date{first version: February 23, 2003; revised: June 16, 2003}

\dedicatory{Dedicated to Mizan Rahman\\
---a pioneer, among other things,\\
of finding new breathtaking identities---\\
on the occasion of his 70th anniversary}

\begin{abstract}
Using a simple classical method we derive bilateral series identities from
terminating ones. In particular, we show how to deduce
Ramanujan's ${}_1\psi_1$ summation from the $q$-Pfaff--Saalsch\"utz
summation. Further, we apply the same method to our previous
$q$-Abel--Rothe summation to obtain, for the first time,
Abel--Rothe type generalizations of Jacobi's triple product identity.
We also give some results for multiple series.
\end{abstract}

\maketitle

\section{Introduction}

Jacobi's~\cite{jacobi} triple product identity,
\begin{equation}\label{jtpi0}
\sum_{k=-\infty}^{\infty}q^{k^2}z^k=
\prod_{j=0}^{\infty}(1-q^{2j+2})(1+zq^{2j+1})(1+z^{-1}q^{2j+1}),
\end{equation}
is one of the most famous and useful identities
connecting number theory and analysis.
Many grand moments in number theory rely on this result, such as
the theorems on sums of squares (cf.\ \cite[Sec.~8.11]{grhyp}),
the Rogers--Ramanujan identities (cf.\ \cite[Sec.~2.7]{grhyp}),
or Euler's pentagonal number theorem (cf.\ \cite[p.~ 51]{bressasm}).
In addition to this identity, different extensions of it,
including Ramanujan's~\cite{hardy} $_1\psi_1$
summation formula (see \eqref{1psi1}) and Bailey's~\cite{bail66}
very-well-poised $_6\psi_6$ summation formula, have 
served as effective tools for studies
in number theory, combinatorics, and physics, see \cite{andappl}.

In this paper, we derive new generalizations of Jacobi's triple
product identity, in one variable and also in severable variables.
Our new extensions look rather unusual. We classify these
to be of ``Abel--Rothe type'', since they are derived
from $q$-Abel--Rothe summations which we previously found in
\cite[Eq.~(8.5)]{krattschloss} and in \cite[Th.~6.1]{schlnammi}.
At the moment, we cannot tell if our new identities have interesting
combinatorial or number-theoretic applications. Nevertheless, we believe
that they are attractive by its own.

Our article is organized as follows. In Section~\ref{secpre}, we
review some basics in $q$-series. In addition of explaining some
standard notation,
we also briefly describe a well-known method employed
in this article for obtaining a bilateral identity from a unilateral
terminating identity, a method already utilized by Cauchy~\cite{cauchy}
in his second proof of Jacobi's triple product identity.
In Section~\ref{sec1psi1}, we apply this classical method to
derive Ramanujan's ${}_1\psi_1$ summation from the
$q$-Pfaff--Saalsch\"utz summation.
According to our knowledge, this very simple proof of the ${}_1\psi_1$
summation has not been given explicitly before. In Section~\ref{secrgj},
we give two Abel--Rothe type generalizations of Jacobi's triple product
identity, see Theorem~\ref{rgj} and Corollary~\ref{rgjc}.
These are consequences of our $q$-Abel--Rothe summation from
\cite[Eq.~(8.5)]{krattschloss}.
In Section~\ref{secmult} we give multidimensional generalizations
of our Abel--Rothe type identities, associated to the root system $A_{r-1}$
(or equivalently, associated to the unitary group $U\!(r)$).
As a direct consequence, we also give an Abel--Rothe type generalization
of the Macdonald identities for the affine root system $A_r$. 
Finally, we establish the conditions of convergence of
our multiple series in Appendix~\ref{secapp}.

{\bf Acknowledgements:}
I would like to thank George Gasper and Mourad Ismail for their comments
on an earlier version of this article. Further, I am grateful to
the anonymous referees for their very detailed comments. In particular,
one of the referees asked me to find Abel--Rothe type extensions of
the Macdonald identities. I thus included such an extension, see
\eqref{armacdid}.
Another referee insisted that I need to clarify the arguments
I had used to prove the convergence of the multiple series in
Theorems~\ref{mrgj} and \ref{mrgjc}. This eventually lead me to
find even slightly more general convergence conditions than
I had originally stated (see Remark~\ref{rem}).

\section{Some basics in $q$-series}\label{secpre}

First, we recall some standard notation for {\em $q$-series} and {\em basic
hypergeometric series}~\cite{grhyp}.
Let $q$ be a (fixed) complex parameter
(called the ``base'') with $0<|q|<1$. Then, for a complex parameter $a$,
we define the {\em $q$-shifted factorial} by
\begin{equation}\label{pochinf}
(a)_{\infty}\equiv(a;q)_{\infty}:=\prod_{j\ge 0}(1-aq^j),
\end{equation}
and 
\begin{equation}\label{poch}
(a)_k\equiv(a;q)_k:=\frac{(a;q)_{\infty}}{(aq^k;q)_{\infty}},
\end{equation}
where $k$ is any integer.
Since we work with the same base $q$ throughout this article,
we can readily omit writing out the base in the $q$-shifted factorials
(writing $(a)_k$ instead of $(a;q)_k$, etc.) as this does not
lead to any confusion.
For brevity, we occasionally employ the condensed notation
\begin{equation*}
(a_1,\ldots,a_m)_k:=(a_1)_k\dots(a_m)_k,
\end{equation*}
where $k$ is an integer or infinity.
Further, we utilize
\begin{equation}\label{defhyp}
{}_s\phi_{s-1}\!\left[\begin{matrix}a_1,a_2,\dots,a_s\\
b_1,b_2,\dots,b_{s-1}\end{matrix}\,;q,z\right]:=
\sum _{k=0} ^{\infty}\frac {(a_1,a_2,\dots,a_s;q)_k}
{(q,b_1,\dots,b_{s-1};q)_k}z^k,
\end{equation}
and
\begin{equation}\label{defhypb}
{}_s\psi_s\!\left[\begin{matrix}a_1,a_2,\dots,a_s\\
b_1,b_2,\dots,b_s\end{matrix}\,;q,z\right]:=
\sum _{k=-\infty} ^{\infty}\frac {(a_1,a_2,\dots,a_s;q)_k}
{(b_1,b_2,\dots,b_s;q)_k}z^k,
\end{equation}
to denote the {\em basic hypergeometric ${}_s\phi_{s-1}$ series},
and the {\em bilateral basic hypergeometric ${}_s\psi_s$ series},
respectively.

A standard reference for basic hypergeometric series
is Gasper and Rahman's text~\cite{grhyp}.
Throughout this article,
in our computations we make decent use of some elementary identities for
$q$-shifted factorials, listed in \cite[Appendix~I]{grhyp}.

We now turn our attention to identities.
One of the simplest summations for basic hypergeometric series
is the terminating $q$-binomial theorem,
\begin{equation}\label{qbin}
{}_1\phi_0\!\left[\begin{matrix}q^{-n}\\
-\end{matrix}\,;q,z\right]=(zq^{-n})_n.
\end{equation}
This can also be written (with $z\mapsto zq^n$) as
\begin{equation}\label{qbino}
\sum_{k=0}^n\begin{bmatrix}n\\k\end{bmatrix}_q(-1)^kq^{\binom k2}z^k
=(z)_n,
\end{equation}
where
\begin{equation}
\begin{bmatrix}n\\k\end{bmatrix}_q:=\frac{(q;q)_n}{(q;q)_k(q;q)_{n-k}}
\end{equation}
denotes the {\em $q$-binomial coefficient}.

Cauchy's~\cite{cauchy} second proof of Jacobi's
triple product identity is very elegant and actually constitutes
a useful method for obtaining bilateral series identities in general.
It is worth looking closely at his proof:
First he replaced in \eqref{qbino} $n$ by $2n$ and then
shifted the summation index $k\mapsto k+n$, which leads to
\begin{equation}
(z)_{2n}=\sum_{k=-n}^n\begin{bmatrix}2n\\n+k\end{bmatrix}_q
(-1)^{n+k}q^{\binom{n+k}2}z^{n+k}.
\end{equation}
Next, he replaced $z$ by $zq^{-n}$ and obtained after some elementary
manipulations
\begin{equation}
(z,q/z)_n=\sum_{k=-n}^n\begin{bmatrix}2n\\n+k\end{bmatrix}_q
(-1)^kq^{\binom k2}z^k.
\end{equation}
Finally, after letting $n\to\infty$ he obtained
\begin{equation}\label{jtpi}
\sum_{k=-\infty}^{\infty}(-1)^kq^{\binom k2}z^k=
(q,z,q/z)_{\infty},
\end{equation}
which is an equivalent form of Jacobi's
triple product identity \eqref{jtpi0}.
 
\section{Ramanujan's ${}_1\psi_1$ summation}\label{sec1psi1}

Hardy~\cite[Eq.~(12.12.2)]{hardy} describes Ramanujan's
${}_1\psi_1$ summation (cf.~\cite[Appendix~(II.29)]{grhyp}),
\begin{equation}\label{1psi1}
_1\psi_1\!\left[\begin{matrix}a\\b\end{matrix}\,;q,z\right]
=\frac{(q,b/a,az,q/az)_{\infty}}{(b,q/a,z,b/az)_{\infty}},
\end{equation}
where $|b/a|<|z|<1$,
as a ``remarkable formula with many parameters''. On the one hand,
it bilaterally extends the nonterminating $q$-binomial theorem
(which is the $b=q$ special case of \eqref{1psi1}), on the other hand
it also contains Jacobi's triple product identity as a special case.
Namely, if in \eqref{1psi1} we replace $z$ by $z/a$, and then let
$a\to\infty$ and $b\to 0$, we immediately obtain \eqref{jtpi}.
Another important special case of \eqref{1psi1} is obtained when $b=aq$,
which is a bilateral $q$-series summation due to Kronecker,
see Weil~\cite[pp.~70--71]{weil}.

Ramanujan (who very rarely gave any proofs)
did not provide a proof for the above summation formula.
It is interesting that Bailey's~\cite[Eq.~(4.7)]{bail66}
very-well-poised ${}_6\psi_6$ summation formula,
although it contains more parameters
than Ramanujan's ${}_1\psi_1$ summation,
does {\em not} include the latter as a special case.
Hahn~\cite[$\kappa=0$ in Eq.~(4.7)]{hahn} independently
established \eqref{1psi1} by considering a first order homogeneous
$q$-difference equation. Hahn thus published the first proof
of the $_1\psi_1$ summation. Not much later, M.~Jackson~\cite[Sec.~4]{mjack}
gave the first elementary proof of \eqref{1psi1}.
Her proof derives the $_1\psi_1$ summation from the $q$-Gau{\ss}
summation, by manipulation of series.
A simple and elegant proof of the $_1\psi_1$ summation formula
was given by Ismail~\cite{ismail}
who showed that the $_1\psi_1$ summation is an immediate consequence of
the $q$-binomial theorem and analytic continuation.

We provide yet another simple proof of the
${}_1\psi_1$ summation formula (which seems to have been unnoticed so far)
by deriving it from the terminating $q$-Pfaff--Saalsch\"utz summation
(cf.~\cite[Eq.~(II.12)]{grhyp}),
\begin{equation}\label{3phi2}
{}_3\phi_2\!\left[\begin{matrix}a,b,q^{-n}\\
c,abq^{1-n}/c\end{matrix}\,;q,q\right]=
\frac{(c/a,c/b)_n}{(c,c/ab)_n}.
\end{equation}

First, in \eqref{3phi2} we replace $n$ by $2n$ and then shift
the summation index by $n$ such that the new sum runs from $-n$ to $n$:
\begin{multline*}
\frac{(c/a,c/b)_{2n}}{(c,c/ab)_{2n}}=
\sum_{k=0}^{2n}\frac{(a,b,q^{-2n})_k}{(q,c,abq^{1-2n}/c)_k}q^k\\
=\frac{(a,b,q^{-2n})_n}{(q,c,abq^{1-2n}/c)_n}q^n\,
\sum_{k=-n}^{n}\frac{(aq^n,bq^n,q^{-n})_k}
{(q^{1+n},cq^n,abq^{1-n}/c)_k}q^k.
\end{multline*}
Next, we replace $a$ by $aq^{-n}$, and we replace $c$ by $cq^{-n}$.
\begin{multline*}
\sum_{k=-n}^{n}\frac{(a,bq^n,q^{-n})_k}
{(q^{1+n},c,abq^{1-n}/c)_k}q^k=
\frac{(c/a,cq^{-n}/b)_{2n}(q,cq^{-n},abq^{1-2n}/c)_n}
{(cq^{-n},c/ab)_{2n}(aq^{-n},b,q^{-2n})_n}q^{-n}\\
=\frac{(c/a)_{2n}(c/b,bq/c,q,q)_n}
{(q)_{2n}(c,q/a,b,c/ab)_n}.
\end{multline*}
Now, we may let $n\to\infty$ (assuming $|c/ab|<1$ and $|b|<1$) while
appealing to Tannery's theorem~\cite{bromwich} for being allowed to
interchange limit and summation. This gives
\begin{equation*}
\sum_{k=-\infty}^{\infty}\frac{(a)_k}
{(c)_k}\left(\frac c{ab}\right)^k=
\frac{(c/a,c/b,bq/c,q)_{\infty}}
{(c,q/a,b,c/ab)_{\infty}},
\end{equation*}
where $|c/a|<|c/ab|<1$.
Finally, replacing $b$ by $c/az$ and then $c$ by $b$ gives \eqref{1psi1}.

\begin{Remark}
The elementary method we use in the above derivation
(exactly the same method already utilized by Cauchy)
has also been exploited by Bailey~\cite[Secs.~3 and 6]{bail66},
\cite{bail22} (see also Slater~\cite[Sec.~6.2]{slater}).
For instance, in \cite{bail22} Bailey applies the method
to Watson's transformation formula of a terminating very-well-poised
$_8\phi_7$ into a multiple of a balanced $_4\phi_3$
\cite[Eq.~(III.18)]{grhyp}. As a result, he obtains a transformation
for a $_2\psi_2$ series, see also Gasper and Rahman~\cite[Ex.~5.11]{grhyp}.
\end{Remark}

\begin{Remark}
We conjecture that any bilateral sum can be obtained from an appropriately
chosen terminating identity by the above method
(as a limit, without using analytic continuation).
However, it is already not known whether Bailey's~\cite[Eq.~(4.7)]{bail66}
$_6\psi_6$ summation formula (cf.~\cite[Eq.~(II.33)]{grhyp})
follows from such an identity.
\end{Remark}

\section{Abel--Rothe type generalizations
of Jacobi's triple product identity}\label{secrgj}

We apply the method of bilateralization\footnote{The Merriam--Webster Online
dictionary ({\tt http://www.m-w.com/cgi-bin/dictionary})
gives for the entry '{\bf -ize}':
``\dots {\bf 1 a} (1){\bf :} cause to be \dots (2){\bf :}
cause to be formed into \dots
{\bf 2 a:} become \dots
{\bf usage}
The suffix {\it -ize} has been productive in English since the time of
Thomas Nashe (1567-1601), who claimed credit for introducing it into
English to remedy the surplus of monosyllabic words.
Almost any noun or adjective can be made into a verb
by adding {\it -ize} $<$hospital{\it ize}$>$ $<$familiar{\it ize}$>$;
many technical terms are coined this way $<$oxid{\it ize}$>$
as well as verbs of ethnic derivation $<$American{\it ize}$>$
and verbs derived from proper names $<$bowdler{\it ize}$>$
$<$mesmer{\it ize}$>$. Nashe noted in 1591 that his coinages in {\it-ize}
were being complained about, and to this day new words in {\it -ize}
$<$final{\it ize}$>$ $<$priorit{\it ize}$>$ are sure to draw
critical fire.''}
we just utilized now to the
following $q$-Abel--Rothe summation \cite[Eq.~(8.5)]{krattschloss},
\begin{equation}\label{qrothegl}
(c)_n=\sum_{k=0}^n\begin{bmatrix}n\\k\end{bmatrix}_q
(1-a-b)\,\big(aq^{1-k}+bq\big)_{k-1}
\big(c(a+bq^k)\big)_{n-k}(-1)^kq^{\binom k2}c^k.
\end{equation}
This summation is different from the $q$-Rothe summation found by
Johnson~\cite[Th.~4]{johnson} which he derived by means of umbral calculus.
It is also different from Jackson's~\cite{jacksonabel} $q$-Abel summation.
Our summation in \eqref{qrothegl} was originally derived in
\cite{krattschloss} by extracting coefficients of a
nonterminating $q$-Abel--Rothe type expansion formula
(actually, the $n\to\infty$ case of \eqref{qrothegl}),
which in turn was derived by inverse relations.
However, it can also be derived directly by inverse relations
(by combining the $q$-Chu--Vandermonde summation with
a specific non-hypergeometric matrix inverse), see
\cite[Sec.~6]{schlnammi}.

In \cite{krattschloss}, \cite{schlnammi}
and \cite{schlnmmi}, we referred to \eqref{qrothegl} as a
$q$-Rothe summation to distinguish it from the $q$-Abel summation
\begin{equation}\label{qabelgl}
1=\sum_{k=0}^n\begin{bmatrix}n\\k\end{bmatrix}_q
(a+b)\big(a+bq^k\big)^{k-1}\big(a+bq^k\big)_{n-k},
\end{equation}
that we derived in \cite[Eq.~(8.1)]{krattschloss}.
It appears that \eqref{qabelgl} is different from any of the
$q$-Abel summations from Jackson~\cite{jacksonabel} or Johnson~\cite{johnson}.
However, it is equivalent to Bhatnagar and
Milne's~\cite{bhatmil} $q$-Abel summation (by reversing the sum).
Above we decided to call \eqref{qrothegl} a $q$-Abel--Rothe summation
since it is also contains \eqref{qabelgl} as a special case.
In fact, if in \eqref{qrothegl} we replace $a$ and $b$ by
$a/c$ and $b/c$, and then let $c\to 0$, we obtain after some algebra
\eqref{qabelgl}.

Our $q$-Abel--Rothe summation in \eqref{qrothegl} is indeed a $q$-extension
of the Rothe summation:
If we divide both sides by $(q)_n$,
do the replacements $a\mapsto q^A-B$, $b\mapsto B$,
$c\mapsto q^{-A-C}$, and then let $q\to 1$, we obtain
Rothe's~\cite{rothe} summation formula
\begin{equation}\label{rothegl}
\binom{A+C}n=\sum_{k=0}^n\frac A{A+Bk}\binom{A+Bk}k\binom{C-Bk}{n-k}.
\end{equation}
Rothe's identity is an elegant generalization of the well-known
Chu--Vandermonde convolution formula, to which it reduces for $B=0$.
Similarly, \eqref{qrothegl} reduces for $b=0$ to the $q$-Chu--Vandermonde
summation listed in Appendix~II, Eq.~(II.6) of \cite{grhyp},
and for $a=0$ to the $q$-Chu--Vandermonde
summation in Appendix~II, Eq.~(II.7) of \cite{grhyp}.

Furthermore, Eq.\ \eqref{qabelgl}
(and thus also the more general \eqref{qrothegl})
is indeed a $q$-extension of the Abel summation:
If in \eqref{qabelgl}
we replace $a$ and $b$, by $\frac A{(A+C)}+\frac B{(A+C)(1-q)}$
and $\frac{-B}{(A+C)(1-q)}$, respectively, and then let $q\to 1$,
we obtain Abel's generalization~\cite{abel} of the $q$-binomial theorem,
\begin{equation}\label{abelgl}
(A+C)^n=\sum_{k=0}^n\binom nk\,A(A+Bk)^{k-1}(C-Bk)^{n-k}.
\end{equation}

Some historical details concerning the Abel and Rothe summations can be found
in Gould~\cite{gouldvc}, \cite{gouldfvc}, and in Strehl~\cite{strehl}.

We now present our main result, an Abel--Rothe type generalization of
\eqref{jtpi}:
\begin{Theorem}\label{rgj}
Let $a$, $b$, and $z$ be indeterminate. Then
\begin{equation}\label{rgjgl}
\frac{(q,z,q/z)_{\infty}}{(1-b)}
=\sum_{k=-\infty}^{\infty}\big(aq^{1-k}+bq\big)_{\infty}
\big(z(a+bq^k)\big)_{\infty}(-1)^kq^{\binom k2}z^k,
\end{equation}
provided $\max(|az|,|b|)<1$.
\end{Theorem}

\begin{proof}
In \eqref{qrothegl}, we first replace $n$ by $2n$, and then shift the
summation index $k\mapsto k+n$. This gives
\begin{multline*}
(c)_{2n}=\sum_{k=-n}^{n}\begin{bmatrix}2n\\n+k\end{bmatrix}_q
(1-a-b)\,\big(aq^{1-n-k}+bq\big)_{k+n-1}\\
\times\big(c(a+bq^{n+k})\big)_{n-k}(-1)^{n+k}q^{\binom{n+k}2}c^{n+k}.
\end{multline*}
Next, we replace $a$ and $c$ by $aq^n$ and $cq^{-n}$. After some
elementary manipulations we obtain
\begin{equation*}
\frac{(c,q/c)_n}{(1-aq^n-b)}=
\sum_{k=-n}^{n}\begin{bmatrix}2n\\n+k\end{bmatrix}_q
\big(aq^{1-k}+bq\big)_{k+n-1}
\big(c(a+bq^k)\big)_{n-k}(-1)^kq^{\binom k2}c^k.
\end{equation*}
Finally, replacing $c$ by $z$ and (assuming $|az|,|b|<1$)
letting $n\to\infty$ while appealing to Tannery's theorem, we formally
arrive at \eqref{rgjgl}. However, it remains to establish the conditions
of convergence. 

Since
\begin{equation*}
\big(aq^{1-k}+bq\big)_{\infty}=(-1)^kq^{-\binom k2}(a+bq^k)^k
\big(1/(a+bq^k)\big)_k\big(aq+bq^{1+k}\big)_{\infty},
\end{equation*}
it is easy to find that if $|az|<1$ then the positive part of the sum,
i.\ e., $\sum_{k\ge 0}$, converges.
Similarly, for the negative part of the sum, i.\ e., $\sum_{k<0}$, we
use
\begin{equation*}
\big(z(a+bq^k)\big)_{\infty}=(-1)^kq^{-\binom k2}z^{-k}\big(aq^{-k}+b\big)^{-k}
\big(q/z(aq^{-k}+b)\big)_{-k}\big(z(aq^{-k}+b)\big)_{\infty},
\end{equation*}
and determine that we need $|b|<1$ for absolute convergence.
\end{proof}

By reversing the sum in \eqref{rgjgl}, we easily deduce the following:

\begin{Corollary}\label{rgjc}
Let $a$, $b$, and $z$ be indeterminate. Then
\begin{equation}\label{rgjcgl}
\frac{(q,zq,1/z)_{\infty}}{(1-az)}
=\sum_{k=-\infty}^{\infty}\big(aq^{-k}+b\big)_{\infty}
\big(zq(a+bq^k)\big)_{\infty}(-1)^kq^{\binom {k+1}2}z^k,
\end{equation}
provided $\max(|az|,|b|)<1$.
\end{Corollary}

\begin{proof}
In \eqref{rgjgl}, we first replace $k$ by $-k$, and then simultaneously
replace $a$, $b$ and $z$, by $bz$, $az$ and $1/z$, respectively.
\end{proof}

We have given analytical convergence conditions for the identities
\eqref{rgjgl} and \eqref{rgjcgl}. However, we would like to point out
that these identities also hold when regarded as identities for
{\em formal power series} over $q$.

It is obvious that both Theorem~\ref{rgj} and Corollary~\ref{rgjc}
reduce to Jacobi's triple product identity \eqref{jtpi} when
$a=b=0$.

If we extract coefficients of $z^n$ on
both sides of \eqref{rgjgl}, using \eqref{jtpi} on the left and
\begin{equation}\label{qexpgl}
(x)_{\infty}=\sum_{j\ge 0}\frac{(-1)^jq^{\binom j2}}{(q)_j}x^j
\end{equation}
(cf.~\cite[Eq.~(II.2)]{grhyp}) on the right hand side, divide
both sides by $(-1)^nq^{\binom n2}$ and replace
$a$ by $aq^{-n}$, we obtain
\begin{equation}\label{extrc}
\frac 1{1-b}=\sum_{j=0}^{\infty}\frac{(aq^j+b)^j}{(q)_j}(aq^{1+j}+bq)_{\infty}
\end{equation}
(which is valid for $|b|<1$), which is, modulo substitution of variables,
our $q$-Abel-type expansion in \cite[Eq.~(3.4)]{schlnammi}.
For a multivariable extension of \eqref{extrc}, see \eqref{mextrc}.
If we now replace $a$ and $b$ by $-BZ$ and $(1-q^A+B)Z$ and then let
$q\to 1^-$, while using $\lim_{q\to 1^-}((1-q)Z)_{\infty}=e^{-Z}$, we recover
Lambert's~\cite{lambert} formula
\begin{equation}\label{lamb}
\frac{e^{AZ}}{1-BZ}=\sum_{j=0}^{\infty}
\frac{(A+Bj)^j}{j!}Z^je^{-BZj},
\end{equation}
which is valid for $|BZe^{1-BZ}|<1$.
Note that in \eqref{lamb}, $Z$ is a redundant parameter. However,
the advantage of writing \eqref{lamb} in this form is that
here we have an identity of power series in the variable $Z$
(having in mind the expansion of the geometric series and of the
exponential function).

In \cite{krattschloss}, \cite{schlnammi}, and
\cite{schlnmmi} we erroneously attributed \eqref{lamb} and
some related expansions to Euler~\cite{euler}, but which are actually
due to Lambert~\cite{lambert}. Nevertheless, Euler's article on
Lambert's identities is significant and is often cited in the literature
as sole reference for these identities
(see e.\ g., P\'olya and Szeg\"o~\cite[pp.~301--302]{polyaszegoe1}).

\section{Multidimensional generalizations}\label{secmult}

Here we extend Theorem~\ref{rgj} and Corollary~\ref{rgjc}
to multiple series associated to the root systems of type A, or
equivalently, associated to the unitary groups. Multiple series,
associated to root systems, or to Lie groups, have been investigated
by various authors. Many different types of such series exist
in the literature. For some results on the special type of series
that are considered in this section, see, e.g., \cite{bhatmil}
\cite{bhatms}, \cite{milne}, \cite{milnschloss}, \cite{rosengrenan},
\cite{schlossmmi}, \cite{schlnammi}, and \cite{schlnmmi}.

In the following, we consider $r$-dimensional
series, where $r$ is a positive integer. For brevity, we
employ the notation $|{\mathbf k}|:=k_1+\dots+k_r$.

If we apply the method of bilateralization to the multidimensional
$q$-Abel--Rothe summations that were derived in \cite{krattschloss}
(see Theorems~8.2 and 8.3 therein), the multiple $q$-Abel--Rothe summations
in Theorems~6.7 and 6.9 of \cite{schlnammi},
or Theorems~3.7 and 3.8 of \cite{schlnmmi},
the resulting series do not converge for higher dimensions. 
The only multidimensional $q$-Abel--Rothe summations we are aware of
that converge when bilateralized are
Theorem~6.11 of \cite{schlnammi}, and (the slightly more general)
Theorem~3.9 of \cite{schlnmmi}. Both these theorems were
derived by applying multidimensional inverse relations,
in particular by combining different higher-dimensional
$q$-Chu--Vandermonde summations with specific multidimensional
non-hypergeometric matrix inverses.

For the sake of simplicity, we consider here only the multilateral identities
arising from Theorem~6.11 of \cite{schlnammi}, a multiple
$q$-Abel--Rothe summation associated to the root system $A_{r-1}$:

\begin{Theorem}\label{rothe3}
Let $a$, $b$, $c$, and $x_1,\dots,x_r$ be indeterminate, and
let $n_1,\dots,n_r$ be nonnegative integers.
Then there holds
\begin{multline}\label{rothe3gl}
(c)_{|{\mathbf n}|}=
\underset{i=1,\dots,r}{\sum_{0\le k_i\le n_i}}\Bigg(
\prod_{i,j=1}^{r}\left[\frac {\left(\frac{x_i}{x_j}q\right)_{n_i}}
{\left(\frac{x_i}{x_j}q\right)_{k_i}
\left(\frac{x_i}{x_j}q^{1+k_i-k_j}\right)_{n_i-k_i}}\right]\\\times
(1-a-b)\,
\big(aq^{1-|{\mathbf k}|}+bq\big)_{|{\mathbf k}|-1}
\big(c(a+bq^{|{\mathbf k}|})\big)_{|{\mathbf n}|-|{\mathbf k}|}
(-1)^{|{\mathbf k}|}q^{\binom{|{\mathbf k}|}2}
c^{|\mathbf k|}.
\Bigg).
\end{multline}
\end{Theorem}

Our multilateral extension of \eqref{rgjgl} is as follows:
\begin{Theorem}\label{mrgj}
Let $a$, $b$, $z$, and $x_1,\dots,x_r$ be indeterminate.
Then
\begin{multline}\label{mrgjgl}
\frac{(z,q/z)_{\infty}}{(1-b)}
\prod_{i,j=1}^{r}\left(\frac{x_i}{x_j}q\right)_{\infty}\\
=\sum _{k_1,\dots,k_r=-\infty}^{\infty}\Bigg(
\prod_{1\le i<j\le r}\left(\frac {x_iq^{k_i}-x_jq^{k_j}}{x_i-x_j}\right)
\big(aq^{1-|{\mathbf k}|}+bq\big)_{\infty}\\\times
\big(z(a+bq^{|{\mathbf k}|})\big)_{\infty}
(-1)^{r|{\mathbf k}|}q^{r\sum_{i=1} ^{r}{\binom{k_i}2}}
z^{|\mathbf k|}\prod_{i=1}^{r}x_i^{rk_i-|{\mathbf k}|}
\Bigg),
\end{multline}
provided $\max(|az|,|b|)<1$.
\end{Theorem}

\begin{proof}
The proof is very similar to the one-dimensional case.
In \eqref{rothe3gl}, we replace $n_i$ by $2n_i$, for $i=1,\dots,r$,
and then shift all the summation indices $k_i\mapsto k_i+n_i$. This gives
\begin{multline*}
(c)_{2|{\mathbf n}|}=
\underset{i=1,\dots,r}{\sum_{-n_i\le k_i\le n_i}}\Bigg(
\prod_{i,j=1}^{r}\left[\frac {\left(\frac{x_i}{x_j}q\right)_{2n_i}}
{\left(\frac{x_i}{x_j}q\right)_{n_i+k_i}
\left(\frac{x_i}{x_j}q^{1+n_i-n_j+k_i-k_j}\right)_{n_i-k_i}}\right]\\\times
(1-a-b)\,\big(aq^{1-|{\mathbf n}|-|{\mathbf k}|}+bq
\big)_{|{\mathbf n}|+|{\mathbf k}|-1}\\\times
\big(c(a+bq^{|{\mathbf n}|+|{\mathbf k}|})\big)_{|{\mathbf n}|-|{\mathbf k}|}
(-1)^{|{\mathbf n}|+|{\mathbf k}|}q^{\binom{|{\mathbf n}|+|{\mathbf k}|}2}
c^{|{\mathbf n}|+|\mathbf k|}\Bigg).
\end{multline*}
We replace $a$, $c$ and $x_i$,
by $aq^{|{\mathbf n}|}$, $cq^{-|{\mathbf n}|}$ and $x_iq^{-n_i}$,
$i=1,\dots,r$, respectively. After some elementary manipulations we obtain
\begin{multline*}
\frac{(c,q/c)_{|{\mathbf n}|}}
{(1-aq^{|{\mathbf n}|}-b)}=
\underset{i=1,\dots,r}{\sum_{-n_i\le k_i\le n_i}}\Bigg(
\prod_{i,j=1}^{r}\left[\frac {\left(\frac{x_i}{x_j}q\right)_{n_i+n_j}}
{\left(\frac{x_i}{x_j}q\right)_{n_j+k_i}
\left(\frac{x_i}{x_j}q^{1+k_i-k_j}\right)_{n_i-k_i}}\right]\\\times
\big(aq^{1-|{\mathbf k}|}+bq\big)_{|{\mathbf n}|+|{\mathbf k}|-1}
\big(c(a+bq^{|{\mathbf k}|})\big)_{|{\mathbf n}|-|{\mathbf k}|}
(-1)^{|{\mathbf k}|}q^{\binom{|{\mathbf k}|}2}
c^{|\mathbf k|}\Bigg).
\end{multline*}
Next, we replace $c$ by $z$ and 
let $n_i\to\infty$, for $i=1,\dots,r$ (assuming
$|az|<1$ and $|b|<1$), while appealing to Tannery's theorem.
Finally, we apply the simple identity
\begin{multline}
\prod_{i,j=1}^{r}\left(\frac{x_i}{x_j}q\right)_{k_i-k_j}
=\prod_{1\le i<j\le r}\left(\frac{x_i}{x_j}q\right)_{k_i-k_j}
\left(\frac{x_j}{x_i}q\right)_{k_j-k_i}\\
=(-1)^{(r-1)|{\mathbf k}|}
q^{-\binom{|{\mathbf k}|}2+r\sum_{i=1}^r\binom{k_i}2}
\prod_{i=1}^rx_i^{rk_i-|{\mathbf k}|}
\prod_{1\le i<j\le r}\left(\frac {x_iq^{k_i}-x_jq^{k_j}}{x_i-x_j}\right)
\end{multline}
and arrive at \eqref{mrgjgl}. For establishing the conditions of
convergence of the series, see Appendix~\ref{secapp}.
\end{proof}

Next, we provide the following multilateral
generalization of \eqref{rgjcgl}.
\begin{Theorem}\label{mrgjc}
Let $a$, $b$, $z$, and $x_1,\dots,x_r$ be indeterminate. Then
\begin{multline}\label{mrgjcgl}
\frac{(zq,1/z)_{\infty}}{(1-az)}
\prod_{i,j=1}^{r}\left(\frac{x_i}{x_j}q\right)_{\infty}\\
=\sum _{k_1,\dots,k_r=-\infty}^{\infty}\Bigg(
\prod_{1\le i<j\le r}\left(\frac {x_iq^{k_i}-x_jq^{k_j}}{x_i-x_j}\right)
\big(aq^{-|{\mathbf k}|}+b\big)_{\infty}\\\times
\big(zq(a+bq^{|{\mathbf k}|})\big)_{\infty}
(-1)^{r|{\mathbf k}|}q^{|{\mathbf k}|+r\sum_{i=1} ^{r}{\binom{k_i}2}}
z^{|\mathbf k|}\prod_{i=1}^{r}x_i^{rk_i-|{\mathbf k}|}
\Bigg),
\end{multline}
provided $\max(|az|,|b|)<1$.
\end{Theorem}

\begin{proof}
In \eqref{mrgjgl}, we first replace $k_i$ by $-k_i$, for $i=1,\dots,r$,
and then simultaneously replace $a$, $b$, $z$ and $x_i$,
by $bz$, $az$, $1/z$ and $1/x_i$, for $i=1,\dots,r$, respectively.
\end{proof}

We have given analytical convergence conditions for the identities
\eqref{mrgjgl} and \eqref{mrgjcgl}. However, as we already observed in
Section~\ref{secrgj} when dealing with the respective one variable cases,
these identities also hold when regarded as identities for
{\em formal power series} over $q$.

We complete this section with an Abel--Rothe type generalization
of the Macdonald identities for the affine root system $A_r$,
as a direct consequence of Theorem~\ref{mrgj}.

If we multiply both sides of \eqref{mrgjgl} by
$\prod_{1\le i<j\le r}(1-x_i/x_j)$, extract the coefficient of $z^M$,
using the Jacobi triple product identity \eqref{jtpi} on the left hand
side and \eqref{qexpgl} on the right hand side, and divide the
resulting identity by $(-1)^Mq^{\binom M2}$, we obtain
\begin{multline}\label{mextrc}
\frac 1{1-b}\,(q)_{\infty}^{r-1}\prod_{1\le i<j\le r}
\left(\frac{x_i}{x_j},\frac{x_j}{x_i}q\right)_{\infty}\\
=\underset{|{\mathbf k}|\le M}{\sum_{k_1,\dots,k_r=-\infty}^{\infty}}
\Bigg(\prod_{1\le i<j\le r}\left(1-\frac{x_i}{x_j}q^{k_i-k_j}\right)
\big(aq^{1-|{\mathbf k}|}+bq\big)_{\infty}\,
\frac{\big(a+bq^{|{\mathbf k}|}\big)^{M-|{\mathbf k}|}}
{(q)_{M-|{\mathbf k}|}}\\\times
(-1)^{(r-1)|{\mathbf k}|}q^{-M|{\mathbf k}|+\binom{|{\mathbf k}|+1}2+
\sum_{i=1}^r r\binom{k_i}2+(i-1)k_i}
\prod_{i=1}^r x_i^{rk_i-|{\mathbf k}|}
\Bigg).
\end{multline}
We use the Vandermonde determinant expansion \eqref{vde} and
a little bit of algebra and obtain
\begin{multline}\label{armacdid}
\frac 1{1-b}\,(q)_{\infty}^{r-1}\prod_{1\le i<j\le r}
\left(\frac{x_i}{x_j},\frac{x_j}{x_i}q\right)_{\infty}\\
=\sum_{\sigma\in{\mathcal S}_r}\operatorname{sgn}(\sigma)
\prod_{i=1}^rx_i^{\sigma(i)-i}
\underset{|{\mathbf k}|\le M}{\sum_{k_1,\dots,k_r=-\infty}^{\infty}}
\Bigg(\big(aq^{1-|{\mathbf k}|}+bq\big)_{\infty}\,
\frac{\big(a+bq^{|{\mathbf k}|}\big)^{M-|{\mathbf k}|}}
{(q)_{M-|{\mathbf k}|}}\\\times
(-1)^{(r-1)|{\mathbf k}|}q^{-M|{\mathbf k}|+\binom{|{\mathbf k}|+1}2+
\sum_{i=1}^r r\binom{k_i}2+(\sigma(i)-1)k_i}
\prod_{i=1}^r x_i^{rk_i-|{\mathbf k}|}
\Bigg).
\end{multline}

If $a=0$ and $b=0$, the terms of the sum in \eqref{armacdid} are
zero unless $|{\mathbf k}|= M$. Specializing this further by setting
$M=0$ gives an identity which has been shown to be equivalent to the
Macdonald identities for the affine root system
$A_r$, see Milne~\cite{milmac}.
Regarding this, we may consider the identity 
\eqref{armacdid} as an Abel--Rothe type generalization of
the Macdonald identities for the affine root system $A_r$.

Concluding, we want to point out that we could have given
an even more general multidimensional Abel--Rothe type
generalization of Jacobi's triple product identity than Theorem~\ref{mrgj},
by multilateralizing Theorem~3.9 of \cite{schlnmmi}
instead of Theorem~6.11 of \cite{schlnammi} as above.
However, we feel that, because of the more complicated factors being
involved, the result would be not as elegant as Theorem~\ref{mrgj}
which is sufficiently illustrative. We therefore decided to refrain
from giving this more general identity.

\begin{appendix}
\section{Convergence of multiple series}\setcounter{equation}{0}
\label{secapp}

Here we prove the conditions of convergence of our multiple series
identity in Theorem~\ref{mrgj} (and thus also of Theorem~\ref{mrgjc}).

We determine the condition for absolute convergence of the
multilateral series in \eqref{mrgjgl} by splitting the entire sum
$\sum_{k_1,\dots,k_r=-\infty}^{\infty}$
into two sums, $\sum_{|{\mathbf k}|\ge 0}$ and
$\sum_{|{\mathbf k}|<0}$, and show the absolute convergence
for each of these separately.

We first consider the sum
$\sum_{|{\mathbf k}|\ge 0}$. We obtain that for $|{\mathbf k}|\ge 0$ the
sum in \eqref{mrgjgl} converges absolutely provided
\begin{multline}\label{bsprothe1}
\underset{|{\mathbf k}|\ge 0}{\sum_{k_1,\dots,k_r=-\infty}^{\infty}}
\bigg|\big(a+bq^{|{\mathbf k}|}\big)^{|{\mathbf k}|}z^{|{\mathbf k}|}
q^{-\binom{|{\mathbf k}|}2+r\sum_{i=1} ^{r}{\binom{k_i}2}}\\\times
\prod_{i=1}^{r}x_i^{rk_i-|{\mathbf k}|}
\prod_{1\le i<j\le r}\big(x_iq^{k_i}-x_jq^{k_j}\big)\bigg|<\infty.
\end{multline}
We use the Vandermonde determinant expansion
\begin{equation}\label{vde}
\prod_{1\le i<j\le r}\big(x_iq^{k_i}-x_jq^{k_j}\big)=
\sum_{\sigma\in{\mathcal S}_r}
\operatorname{sgn}(\sigma)\prod_{i=1}^rx_i^{r-\sigma(i)}
q^{(r-\sigma(i))k_i},
\end{equation}
where ${\mathcal S}_r$ denotes the symmetric group of order $r$,
interchange summations in \eqref{bsprothe1} and obtain
$r!$ multiple sums each corresponding to a permutation
$\sigma\in{\mathcal S}_r$. Thus for $|{\mathbf k}|\ge 0$
the series in \eqref{mrgjgl}
converges provided
\begin{equation}\label{bsprothe2}
\underset{|{\mathbf k}|\ge 0}{\sum_{k_1,\dots,k_r=-\infty}^{\infty}}
\left|\big(a+bq^{|{\mathbf k}|}\big)^{|{\mathbf k}|}z^{|{\mathbf k}|}
q^{-\binom{|{\mathbf k}|}2+r\sum_{i=1}^r\binom{k_i}2}
\prod_{i=1}^{r}q^{(r-\sigma(i))k_i}x_i^{rk_i-|{\mathbf k}|}
\right|<\infty,
\end{equation}
for any $\sigma\in{\mathcal S}_r$.

The next step is crucial and typically applies in the theory of
multidimensional basic hypergeometric series over the root system
$A_{r-1}$ for a class of series.
(For instance, it applies to several of the multilateral
summations and transformations in \cite{milnschloss}.)
In the summand of \eqref{bsprothe2}, we have
\begin{equation}\label{sid}
-\binom{|{\mathbf k}|}2+r\sum_{i=1}^r\binom{k_i}2=
-\frac {(r-1)}2|{\mathbf k}|+\frac 12\sum_{1\le i<j\le r}(k_i-k_j)^2
\end{equation}
appearing in the exponent of $q$. Since $|{\mathbf k}|\ge 0$,
we can assume, without loss of generality, that in
particular $k_r\ge 0$. (At least one summation index is nonnegative
-- choose it to be the $r$-th, by relabelling, if necessary.)

In order to exploit the quadratic powers of $q$ in the sum
(which contribute to convergence), we make the substitutions
\begin{equation*}
k_i\mapsto\sum_{i\le l\le r}m_l,\qquad\mbox{for $i=1,\dots,r$}.
\end{equation*}
Under these substitutions $|{\mathbf k}|$ becomes $\sum_{l=1}^rlm_l$,
while for $i<j$, $k_i-k_j$ becomes $\sum_{i\le l<j}m_l$.
We now use \eqref{sid} and replace
$q^{-\binom{|{\mathbf k}|}2+r\sum_{i=1}^r\binom{k_i}2}$ by
\begin{equation*}
q^{-\frac {(r-1)}2\sum_{l=1}^rlm_l+\frac 12\sum_{l=1}^{r-1}m_l^2}
\end{equation*}
(we left out some quadratic powers),
for comparison with a dominating multiple series, and obtain that 
for $|{\mathbf k}|\ge 0$ and $k_r\ge 0$ the series in
\eqref{mrgjgl} converges provided
\begin{multline*}
\underset{m_r,\sum_{l=1}^rlm_l\ge 0}
{\sum_{m_1,\dots,m_r=-\infty}^{\infty}}
\bigg|\big(a+bq^{\sum_{l=1}^rlm_l}\big)^{\sum_{l=1}^rlm_l}
z^{\sum_{l=1}^rlm_l}
q^{-\frac {(r-1)}2\sum_{l=1}^rlm_l+\frac 12\sum_{l=1}^{r-1}m_l^2}\\\times
\prod_{i=1}^{r}q^{(r-\sigma(i))\sum_{i\le l\le r}m_l}
x_i^{r\sum_{i\le l\le r}m_l-\sum_{1\le l\le r}lm_l}
\bigg|<\infty,
\end{multline*}
for any $\sigma\in{\mathcal S}_r$.
The above series converges absolutely if
\begin{multline}\label{bsprothe3}
\underset{m_r,\sum_{l=1}^rlm_l\ge 0}
{\sum_{m_1,\dots,m_r=-\infty}^{\infty}}
\bigg|\big(azq^{-\frac {(r-1)}2}{\textstyle\prod_{i=1}^{r}
x_i^{-1}}\big)^{\sum_{l=1}^rlm_l}
q^{\frac 12\sum_{l=1}^{r-1}m_l^2}\\\times
\prod_{i=1}^{r}q^{(r-\sigma(i))\sum_{i\le l\le r}m_l}
x_i^{r\sum_{i\le l\le r}m_l}
\bigg|<\infty.
\end{multline}
Now, the series in \eqref{bsprothe3} is dominated by
\begin{multline}
\prod_{l=1}^{r-1}\sum_{m_l=-\infty}^{\infty}
\bigg|\big(azq^{-\frac{r-1}2}{\textstyle\prod_{i=1}^rx_i^{-1}}\big)^{lm_l}
q^{\frac 12 m_l^2}\,q^{\sum_{1\le i\le l}(r-\sigma(i))m_l}
\prod_{1\le i\le l}x_i^{rm_l}\bigg|\\\times
\sum_{m_r=0}^{\infty}
\bigg|\big(azq^{-\frac{r-1}2}{\textstyle\prod_{i=1}^rx_i^{-1}}\big)^{rm_r}
q^{\sum_{1\le i\le r}(r-\sigma(i))m_r}
\prod_{1\le i\le r}x_i^{rm_r}\bigg|.
\end{multline}
We deduce by d'Alembert's ratio test that the product of the first
$r-1$ series converge everywhere due to the
quadratic powers of $q$ (since $|q|<1$). 
Further, by the same test, the $r$-th series converges whenever
\begin{equation*}
\bigg|\big(azq^{-\frac{r-1}2}
{\textstyle\prod_{i=1}^rx_i^{-1}}\big)^r
q^{\sum_{1\le i\le r}(r-\sigma(i))}{\textstyle\prod_{i=1}^rx_i^r}\bigg|
=|az|^r<1,
\end{equation*}
or equivalently, whenever $|az|<1$.
 
The absolute convergence of the sum $\sum_{|{\mathbf k}|< 0}$ is
established in a similar manner. In this case we obtain that for
$|{\mathbf k}|< 0$ the sum in \eqref{mrgjgl} converges absolutely provided
\begin{multline}\label{bsprothe4}
\underset{|{\mathbf k}|<0}{\sum_{k_1,\dots,k_r=-\infty}^{\infty}}
\bigg|\big(aq^{-|{\mathbf k}|}+b\big)^{-|{\mathbf k}|}
q^{-\binom{|{\mathbf k}|}2+r\sum_{i=1} ^{r}{\binom{k_i}2}}\\\times
\prod_{i=1}^{r}x_i^{rk_i-|{\mathbf k}|}
\prod_{1\le i<j\le r}\big(x_iq^{k_i}-x_jq^{k_j}\big)\bigg|<\infty.
\end{multline}
The further analysis is as follows. We use \eqref{vde} and \eqref{sid}
and assume that for $|{\mathbf k}|<0$, without loss of generality,
$k_r<0$.
In a very similar analysis to above one easily finds the condition
$|b|<1$ for absolute convergence.

\begin{Remark}\label{rem}
In an earlier version of this article we had given a smaller region of
convergence for the series in \eqref{mrgjgl}. In particular,
instead of
\begin{equation*}
\max(|az|,|b|)<1
\end{equation*}
we had given the condition
\begin{equation}\label{oldrel}
|az|<\big|q^{\frac{r-1}2}x_j^{-r}
{\textstyle\prod_{i=1}^rx_i}\big|<|q^{r-1}b^{-1}|,
\qquad\mbox{for $j=1,\dots,r$}.
\end{equation}
To see that this gives a smaller region of convergence
(for $r>1$), assume we would have instead $\max(|az|,|b|)\ge 1$.
Now take the ``product'' of the whole relation \eqref{oldrel}
over all $j=1,\dots,r$. This gives
$|az|^r<|q|^{\binom r 2}<|q^{r(r-1)}b^{-r}|$, or equivalently,
after taking $r$-th roots,
$|az|<|q|^{\frac{r-1}2}<|q^{r-1}b^{-1}|$, which apparently
contradicts $\max(|az|,|b|)\ge 1$ since $r\ge 1$.

The same argument which leads to the convergence condition
in Theorem~\ref{mrgj} (in contrast to the condition \eqref{oldrel})
can be used to improve some analogous results
given in the literature. In particular, this concerns
the papers \cite{milne}, \cite{milnschloss}, \cite{schlnammi},
\cite{schlnmmi} (and possibly others).
\end{Remark}

\end{appendix}

\end{document}